\documentclass[a4paper,12pt]{article}
\usepackage{amsmath}
\usepackage{amssymb}
\usepackage{amsthm}
\usepackage{mathtools}
\usepackage{suffix}
\usepackage{stmaryrd}
\usepackage{bbold}
\usepackage{breqn}
\usepackage{txfonts}
\usepackage{graphicx}
\usepackage{xfrac}
\usepackage{booktabs}

\usepackage[dvipsnames]{xcolor}
\usepackage{soul}
\usepackage{comment}
\specialcomment{coloredcomment}{\begingroup\color{blue!50}}{\endgroup}

\usepackage{accents}

\usepackage[english]{babel}
\usepackage{blindtext}

\DeclarePairedDelimiter\floor{\lfloor}{\rfloor}

\usepackage{bm}
\bmdefine{\ba}{a}
\bmdefine{\bb}{b}
\bmdefine{\bc}{c}
\bmdefine{\bd}{d}
\bmdefine{\be}{e}
\bmdefine{\bf}{f}
\bmdefine{\bg}{g}
\bmdefine{\bj}{j}
\bmdefine{\bk}{k}
\bmdefine{\bn}{n}
\bmdefine{\bp}{p}
\bmdefine{\bq}{q}
\bmdefine{\br}{r}
\bmdefine{\bx}{x}
\bmdefine{\by}{y}
\bmdefine{\bu}{u}
\bmdefine{\bv}{v}
\bmdefine{\bw}{w}
\bmdefine{\bz}{z}
\bmdefine{\bE}{E}
\bmdefine{\bB}{B}
\bmdefine{\bC}{C}
\bmdefine{\bD}{D}
\bmdefine{\bJ}{J}
\bmdefine{\bR}{R}
\bmdefine{\bS}{S}
\bmdefine{\bT}{T}
\bmdefine{\bvarphi}{\varphi}
\bmdefine{\bphi}{\phi}
\bmdefine{\bpsi}{\psi}
\bmdefine{\bchi}{\chi}
\bmdefine{\bPsi}{\Psi}
\bmdefine{\blambda}{\lambda}
\bmdefine{\bmu}{\mu}
\bmdefine{\bnu}{\nu}

\def\xfoo#1^#2\relax#3\valign{%
\mathbf{#1}\ifx\valign#2\valign\else^{\mathbf{#2}}\fi}

\usepackage[mathscr]{euscript}

\makeatletter
\newcommand*{\defeq}{\mathrel{\rlap{%
					 \raisebox{0.3ex}{$\m@th\cdot$}}%
					 \raisebox{-0.3ex}{$\m@th\cdot$}}%
					 =}
\makeatother

\makeatletter
\newcommand*{\eqdef}{=\mathrel{\rlap{%
					 \raisebox{0.3ex}{$\m@th\cdot$}}%
					 \raisebox{-0.3ex}{$\m@th\cdot$}}%
					 }
\makeatother

\newcommand{\op}[1]{\operatorname{#1}}
\newcommand{\uo}{\operatorname{\mathfrak{u}}}
\newcommand{\bbb}[1]{\mathbb{#1}}
\newcommand{\llb}{\llbracket}
\newcommand{\rrb}{\rrbracket}

\def\XXint#1#2#3{{\setbox0=\hbox{$#1{#2#3}{\int}$}
     \vcenter{\hbox{$#2#3$}}\kern-.5\wd0}}

\newcommand\wrapped[1]%
  {%
   \begin{array}{@{}l@{}}#1\end{array}%
  }
  
\theoremstyle{definition}

\theoremstyle{definition}

\theoremstyle{definition}

\DeclarePairedDelimiterX\braket[2]{\langle}{\rangle}{#1\,\delimsize\vert\,\mathopen{}#2}

\DeclarePairedDelimiterX\MeijerM[3]{\lparen}{\rparen}%
{#3\, \delimsize\vert\,\begin{smallmatrix}#1 \\ #2\end{smallmatrix}}

\newcommand\MeijerG[8][]{%
  \mathbf{G}^{\,#2,#3}_{\,#4,\,#5}\MeijerM[#1]{#6}{#7}{#8}}

\WithSuffix\newcommand\MeijerG*[7]{%
  \mathbf{G}^{\,#1,#2}_{\,#3,\,#4}\MeijerM*{#5}{#6}{#7}}

\DeclarePairedDelimiterX\FoxM[3]{\lbrack}{\rbrack}%
{#3\, \delimsize\vert\begin{smallmatrix}#1 \\ #2\end{smallmatrix}}

\newcommand\FoxH[8][]{%
  \mathbf{H}^{\,#2,#3}_{\,#4,\,#5}\FoxM[#1]{#6}{#7}{#8}}

\WithSuffix\newcommand\FoxH*[7]{%
  \mathbf{H}^{\,#1,#2}_{\,#3,\,#4}\FoxM*{#5}{#6}{#7}}

\usepackage[]{hyperref}
\hypersetup{
    bookmarks=true,
    bookmarksnumbered=true,     
    bookmarksopen=true,         
    colorlinks=true,            
    pdfstartview=Fit,           
    pdfpagemode=UseOutlines,    
    pdfpagelayout=TwoPageRight
}

\usepackage{cleveref}
\crefformat{pluralequation}{#2eqs.~(#1)#3}
\Crefformat{pluralequation}{#2Eqs.~(#1)#3}

\usepackage{authblk}

\providecommand{\keywords}[1]
{
  \small	
  \textbf{\textit{Keywords---}} #1
}

\begin{document}
	
\title{Le Roy, Lerch and Legendre chi functions and generalised Borel-Le Roy transform}
\author{Giuseppe Dattoli\footnote{E-mail: pinodattoli@libero.it (Giuseppe Dattoli)}}
\author{Roberto Ricci\footnote{E-mail: roberto.ricci@enea.it (Roberto Ricci)}}
\affil{ENEA, Nuclear Department, Frascati Research Center, Via E. Fermi 45, 00044 Frascati (Rome), Italy}
\date{}
\setcounter{Maxaffil}{0}
\renewcommand\Affilfont{\itshape\small}
\maketitle

\section*{Abstract}

The Le Roy function has been the focus of intensive research in recent years, owing both to its relevance in analysis and its versatility in applications involving fractional differential operators. Other special functions -- such as the Lerch transcendent and the Legendre chi function -- have found applications ranging from Bose-Einstein and Fermi-Dirac statistics in physics to pure mathematical investigations involving polylogarithms and Dirichlet L-series.

In this article, we present a unified framework based on a recent reformulation of Indicial Umbral Theory (IUT) grounded in the formal theory of power series. Within this setting, we study the properties and generalisations of these special functions. In particular, we build upon the revised formulation of IUT to incorporate the role of the Borel-Le Roy transform, and to explore the extension of the formalism to divergent series via appropriate resummation techniques.
\vskip0.2cm
\keywords{indicial umbral theory, special functions, Borel-Le Roy transform}

\section*{Introduction}\label{sec:introduction}

The theory of special functions has undergone a series of evolutionary stages. It initially arose from the study of the analytic properties of certain functions and polynomials deemed "special", and gradually developed into a pursuit of a unified conceptual framework. Authoritative textbooks offer a synthesis of these developments (for a partial list, see refs. \cite{CourantHilbert1953, MorseFeshbach1953, InfeldHull1951, Rainville1960, Lebedev1965, Miller1964}), which emerged within specific contexts shaped by the prevailing mathematical paradigms of their time.

A central challenge in this evolution has been identifying the property that qualifies a function as "special". A pivotal moment in this discourse was marked by Wigner's perspective, articulated in his Princeton lectures, which catalysed a paradigm shift: the reinterpretation of special functions as matrix elements of Lie group representations \cite{Wigner1955, Vilenkin1965}. This group-theoretic formulation opened a new path for unification.

The work of Talman \cite{Talman1968} represented a significant advance in this direction, offering a seminal contribution that paved the way for further investigations \cite{Gilmore2006, Wasson2013}.

It became increasingly clear that special functions are specific solutions to families of ordinary differential equations (ODEs) with non-constant coefficients \cite{Andrews1999}. This realisation provided the foundation for alternative unifying strategies. In particular, the hypergeometric differential equation and its associated solutions have emerged as especially well-suited to this goal.

An additional unifying framework has been offered by the Umbral Calculus (UC), whose conceptual roots span nearly two centuries \cite{JordanBell1965}. The theory originated from the observation of formal analogies among seemingly unrelated special polynomials, which -- under suitable conditions -- could be treated as ordinary monomials \cite{Dowker2013}. This insight culminated in Steffensen's introduction of poweroids \cite{Steffensen1941}, and later in the development of the Quasi-Monomiality formalism by Dattoli and Torre \cite{Dattoli2000}.

In the 1970s, Roman and Rota provided a rigorous foundation for the Umbral Calculus by applying the theory of linear functionals. They demonstrated that UC could be formalised as an algebra generated by linear functionals on the vector space of polynomials in a variable $z$ \cite{RomanRota1973}.

More recently, UC has evolved into the Indicial Umbral Theory (IUT) \cite{Babusci2019, ricci2026rigorous}, which incorporates special functions and introduces effective computational tools through a redefinition of umbral techniques and a refined use of Borel-Laplace transforms \cite{ricci2026rigorous, DattoliLicciardi2020}. In particular, ref. \cite{ricci2026rigorous} presents a reformulation of the original IUT (initially termed Indicial Umbral Calculus, IUC), grounded in the differential algebra of formal power series \cite{MarinoResurgence}. This framework enables a mathematically sound definition of the umbral operator -- a perspective reminiscent of Roman and Rota's algebraic approach, but with a significantly greater emphasis on analytic properties.

A notable byproduct of this formulation is its ability to extend the umbral formalism to divergent formal power series, interpreting them as asymptotic expansions that can be resummed via generalised versions of the Borel-Laplace transform.

In this article, we apply IUT techniques to the study of the Le Roy, Lerch, and Legendre functions. We demonstrate that these functions serve as valuable benchmarks for testing the potential of the method and uncovering new structural properties of the functions themselves.

\section{Le Roy function}\label{sec:umbral_le_roy}

The function defined by the series
\begin{dmath}\label{eq:leroy}
	L(\zeta; \mu) \defeq \sum_{r=0}^\infty \frac{\zeta^r}{\Gamma(1+r)^{\,\mu}} \condition*{\mu \in \bbb{R}}
\end{dmath}.
was introduced by Le Roy at the beginning of the last century \cite{LeRoy1899} and  employed in studies regarding the asymptotic properties of the analytic continuation of the sum of power series.

More recently, several authoritative articles (for a partial list see \cite{Gerhold2012, GarraPolito2013, Garrappa2017, Kiryakova2022, PanevaKonovska2023, MeherezRaza2025} and references therein) have renewed the interest in this function, by establishing generalisations, opening new fields of research -- e.g. in fractional calculus -- and proposing novel applications. As a significative example, in a paper dedicated to the study of stochastic differential equations, Kolokoltsov \cite{Kolokoltsov2021} considered the function $L(\zeta;1/2)$, stating that this function plays the same role for stochastic equations as the ordinary exponential for the deterministic case.

In order to investigate the properties of the Le Roy function in the context of the indicial umbral theory (IUT) \cite{Babusci2019, LicciardiDattoli2022}, we introduce the following class of ground states:
\begin{dmath}\label{eq:phi_general}
	\phi_{\alpha,\, \beta}^{\,\mu}(t) \defeq {
	\frac{1}{\Gamma(\beta+\alpha t)^{\,\mu}} 
	\condition*{\alpha >0,\; \beta \in \bbb{C},\;\mu \in \bbb{R} }
}
\end{dmath},
which is an obvious generalisation of the commonly used class:
\begin{dmath}
	\phi_{\alpha,\, \beta}(t) \defeq {
	\frac{1}{\Gamma(\beta+\alpha t)} 
	\condition*{\alpha >0,\; \beta \in \bbb{C},\;\mu \in \bbb{R} }
}
\end{dmath}.
We also consider the specialisation of \cref{eq:phi_general} obtained by setting $\alpha = \beta = 1$, i.e. the particular class of ground states:
\begin{dmath}\label{eq:phi_mu}
	\phi^{\,\mu}(t) \defeq {
	\phi_{1,1}^{\,\mu}(t) =
	\frac{1}{\Gamma(1 + t)^{\,\mu}} 
	\condition*{\mu \in \bbb{R} }
}
\end{dmath}.

By using \cref{eq:phi_mu}, it is possible to express \cref{eq:leroy} in the following simple umbral form:
\begin{dmath}\label{eq:leroy_umbral}
	L(\zeta; \mu) = {
		\frac{1}{1 - \zeta\uo}[\phi^{\,\mu}] =
		\op{e}^{\zeta\uo}\,[\phi^{\,\mu -1}]
	}
\end{dmath}.

On the other hand, by applying the same umbral operator to the general ground state \cref{eq:phi_general}, we obtain one of the possible generalisations of the Le Roy function \cite{CurcioDattoliDiPalma2025}. Namely:
\begin{dmath}\label{eq:generalised_leroy_umbral}
	L(\zeta; \alpha, \beta, \mu) = {
		\op{e}^{\zeta\uo}\,[\phi_{\alpha,\, \beta}^{\,\mu-1}] \defeq
		\sum_{r=0}^\infty \frac{\zeta^r}{r!\Gamma(1 + \beta + \alpha r)^{\,\mu-1}}
	}
\end{dmath}.
Both umbral identity \cref{eq:leroy_umbral} and \cref{eq:generalised_leroy_umbral} are easily proved by expanding the exponential in Maclaurin series and exploiting the definition of $\uo$ as a functional in the differential subalgebra of analytically converging formal series \cite{ricci2026rigorous}, namely:
\begin{dmath}
	\uo^r [\varphi] \defeq {
		\varphi(r) \condition*{\varphi \in \bbb{C}\{ t\} \subset \bbb{C}\llb t \rrb}
	}
\end{dmath},
provided $r$ is within the domain of convergence of $\varphi$.
For further comments see \cite{Babusci2019, LicciardiDattoli2022, CurcioDattoliDiPalma2025, Dattoli2017}.

The use of the previous umbral restyling enables a significant simplification of the study of the properties of Le Roy function and its generalisations. Taking e.g. the repeated derivative with respect to $\zeta$ of both sides of \cref{eq:leroy_umbral}, one eventually finds:
\begin{dmath}\label{eq:leroy_umbral_n_derivative}
	\partial^n_\zeta\,L(\zeta; \mu) = {
		\uo^n \op{e}^{\zeta\uo}\,[\phi^{\,\mu -1}] =
		\sum_{r=0}^\infty \frac{\zeta^r}{r! \Gamma(1 + n + r)^{\,\mu -1}} =
		L(\zeta; 1, n, \mu)
	}
\end{dmath},
where use has been made of the definition \labelcref{eq:generalised_leroy_umbral}.

IUT methods also enable the evaluation of infinite integrals involving the Kolokoltsov function (for further comments see \cite{Babusci2019, LicciardiDattoli2022}), for example:
\begin{dmath}
	\int_{-\infty}^\infty	\mathrm{d}\zeta\, L(-\zeta^2; \mbox{\small $\frac{1}{2}$}) = {
		\int_{-\infty}^\infty	\mathrm{d}\zeta\, \op{e}^{-\zeta^2\uo}\,[\phi^{-1/2}] =
		\sqrt{\pi}\uo ^{-1/2}[\phi^{-1/2}] =
		\sqrt[4]{\pi^3}
	}.
\end{dmath}

The same formalism can be employed to study the three-parameter Mittag-Leffler a.k.a. Prabhakar function \cite{Prabhakar1971}:
\begin{dmath}\label{eq:prabhakar}
	E_{\alpha,\,\beta,\,\gamma}(\zeta) = {
		\sum_{r = 0}^\infty \frac{(\gamma)_r}{r!\,\Gamma(\alpha r + \beta)}\, \zeta^r \condition*{\alpha >0,\; \beta, \gamma \in \bbb{C}}
	}
\end{dmath}.
Using the ground state 
\begin{dmath*}
		\psi_{\alpha,\,\beta\,,\gamma}(t) = {
			\frac{(\gamma)_t}{\Gamma(\alpha t + \beta)} \equiv
			\frac{\Gamma(\gamma + t)}{\Gamma(\gamma)\,\Gamma(\alpha t + \beta)}
			\condition*{\alpha >0,\; \beta, \gamma \in \bbb{C}}
		}
		\end{dmath*},
		where $(\gamma)_t \defeq \Gamma(\gamma + t)/\Gamma(\gamma)$ is a generalised version of the rising Pochhammer symbol, \cref{eq:prabhakar} can be expressed in the umbral form:
\begin{dmath}\label{eq:umbral_prabhakar}
	E_{\alpha,\,\beta,\,\gamma}(\zeta) =  {
		\op{e}^{\zeta\uo} \,[\psi_{\alpha,\, \beta,\,\gamma}]  
	}
\end{dmath}.

An example of the utility of the IUT formalism is provided by the evaluation of the successive derivatives of \cref{eq:prabhakar}. Applying the $\partial^n_\zeta$ operator to both sides of \cref{eq:umbral_prabhakar} and proceeding as before we obtain:
\begin{dmath}\label{eq:umbral_prabhakar_n_derivative}
	\partial^n_\zeta\,E_{\alpha,\,\beta,\,\gamma}(\zeta) =  {
		\uo^n \op{e}^{\zeta\uo} \,[\psi_{\alpha,\, \beta,\,\gamma}] =
		\sum_{r = 0}^\infty \frac{\zeta^r\uo^{r+n}}{r!} \,[\psi_{\alpha,\, \beta,\,\gamma}] 
	} = {
		\sum_{r = 0}^\infty \frac{\Gamma(\gamma + r + n)}{\Gamma(\gamma)\,\Gamma(\alpha r + \beta + \alpha n)}\, \frac{\zeta^r}{r!} =
		(\gamma)_n E_{\alpha,\,\beta +\alpha n,\,\gamma + n}(\zeta)
	}
\end{dmath}.

Another example of application involving integration is (see \cite{GonzalezMoll2016} for a discussion on Pochhammer symbols with negative index):
\begin{dmath}
	\int_{-\infty}^\infty	\mathrm{d}\zeta\, E_{\alpha,\,\beta,\,\gamma}(-\zeta^2) = {
		\int_{-\infty}^\infty	\mathrm{d}\zeta\, \op{e}^{-\zeta^2\uo}\,[\psi_{\alpha,\, \beta,\,\gamma}]
		} = {
		\sqrt{\pi}\uo ^{-1/2}[\psi_{\alpha,\, \beta,\,\gamma}] =
		\sqrt{\pi}\frac{(\gamma)_{-1/2}}{\Gamma(\beta -\alpha/2)}
	}
\end{dmath}.

Before concluding this opening section, we introduce a further element of discussion, whose role is emphasised in the final part of the article.

Let us consider the \emph{integral Borel transform} of the Leroy function, defined as:
\begin{dmath}\label{eq:integral_borel_transform}
	L^{(\mu)}_{\mathcal{B}}(\zeta) \equiv 
	{
		\mathcal{B}[L^{(\mu)}](\zeta) \defeq 
	\int_0^\infty	\mathrm{d}t\, \op{e}^{-t} L^{(\mu)}(\zeta t) 
	}
\end{dmath},
where we have adopted for convenience the alternative notation $L^{(\mu)}(\zeta) \equiv L(\zeta; \mu)$.
Note that 
\begin{dmath}
	\int_0^\infty	\mathrm{d}t\, \op{e}^{-t} L^{(\mu)}(\zeta t) = {
	\frac{1}{\zeta}\int_0^\infty	\mathrm{d}x\, \op{e}^{-x/\zeta} L^{(\mu)}(x) \equiv
	\frac{1}{\zeta}\,\mathcal{L}_1[L^{(\mu)}](\zeta)
	}
\end{dmath},
i.e. $\zeta \,L^{(\mu)}_{\mathcal{B}}(\zeta)$ coincides with the generalised Laplace transform of order 1 of the function $L^{(\mu)}$. On the other hand, as a power series, 
\begin{dmath}
	L^{(\mu)}(\zeta) = \op{B}_1[\tilde{L}^{(\mu-1)}](\zeta)
\end{dmath},
where $\tilde{L}^{(\mu)}(\zeta) \defeq \zeta\,L^{(\mu)}(\zeta)$ and $\op{B}_1$ is the \emph{formal Borel transform} operator of order 1 \cite{ricci2026rigorous}. Since the series $\tilde{L}^{(\mu-1)}$ converges to an entire function, we know from the theory of Borel-Laplace resummation that $\mathcal{L}_1[\op{B}_1[\tilde{L}^{(\mu-1)}]](\zeta) = \tilde{L}^{(\mu-1)}(\zeta)$. It follows that:
 \begin{dmath}\label{eq:leroy_identity}
	L^{(\mu)}_{\mathcal{B}}(\zeta)  = L^{(\mu-1)}(\zeta)
\end{dmath}.

\noindent 
\Cref{eq:leroy_identity} can be easily verified by directly solving the integral in \cref{eq:integral_borel_transform}.

A similar argument shows that the following Borel-Le Roy transform of the generalised Le Roy function \cref{eq:generalised_leroy_umbral}:
\begin{dmath}
	L_{\op{BL}}(\zeta; \alpha, \beta, \mu) \defeq 
	{
	\int_0^\infty	\mathrm{d}t\, \op{e}^{-t} t^{\,\mu} L(\zeta t^{\,\beta}; \alpha, \beta, \mu) 
	}
\end{dmath},
satisfies the identity:
\begin{dmath}
	L_{\op{BL}}(\zeta; \alpha, \beta, \mu) = L(\zeta; \alpha, \beta, \mu -1)
\end{dmath}.
Further examples will be discussed in the section devoted to final comments.

In this introductory section we have touched on the formalism we will use in the forthcoming part of the article, where we frame in the IUT context the Lerch transcendent and Legendre chi function.

\section{Lerch transcendent}\label{sec:umbral_lerch}

The Lerch transcendent \cite{Lerch1887} is defined by the series:
\begin{dmath}\label{eq:lerch}
	\Phi(\zeta; \alpha, s) \defeq \sum_{r=0}^\infty \frac{\zeta^r}{(r + \alpha)^s}
\end{dmath},
converging for any $\alpha > 0$ in $|\zeta|\; < 1$, and also in $|\zeta|\; = 1$ if $\op{Re} s > 0.$

The reason for the interest in this function stems from its association with the polylogarithm function \cite{Koelbig1996}, Dirichlet $\eta$ and Riemann-Hurwitz $\zeta$ functions \cite{ZetaDirichletGroupRef}.
By introducing the new classes of ground states
\begin{dgroup}
	\begin{dmath}\label{eq:generalised_lerch_ground_state}
		\nu_{\alpha,  \,\beta, s}(t) \defeq 
			\frac{(\beta)_t}{(t + \alpha)^s}
	\end{dmath},
	\begin{dmath}\label{eq:lerch_ground_state}
		\nu_{\alpha, s}(t) \defeq {
			\nu_{\alpha, 1, s}(t) =
			\frac{\Gamma(1+t)}{(t + \alpha)^s}
		}
	\end{dmath},
\end{dgroup}
it is evident that \cref{eq:lerch} can be written in umbral form as:
\begin{dmath}\label{eq:umbral_lerch}
	\Phi(\zeta; \alpha, s) = \op{e}^{\zeta\uo} \,[\nu_{\alpha, s}]
\end{dmath}.

By applying the same umbral operator to the general ground state \cref{eq:generalised_lerch_ground_state}, we obtain the generalised Lerch transcendent:
\begin{dmath}\label{eq:generalised_lerch}
	\Phi(\zeta; \alpha, \beta, s) \defeq {
		\sum_{r=0}^\infty \frac{(\beta)_r}{(r + \alpha)^s}\,\frac{\zeta^r}{r!} =
		\op{e}^{\zeta\uo} \,[\nu_{\alpha, \,\beta, s}] 
	}
\end{dmath}, 
which reduces to \cref{eq:umbral_lerch} for $\beta = 1$. 

If we are interested to the properties under derivative of these functions, we can take again advantage from their exponential umbral images and find:

\begin{dgroup}
	\begin{dmath}\label{eq:umbral_lerch_n_derivative}
		\partial^n_\zeta \, \Phi(\zeta; \alpha, s) = {
		\uo^n \op{e}^{\zeta\uo} \,[\nu_{\alpha, s}] = 
		\sum_{r=0}^\infty \frac{\Gamma(1+r+n)}{(r + n + \alpha)^s}\,\frac{\zeta^r}{r!}
	}
		= (1)_n\,\Phi(\zeta; \alpha + n, 1 + n, s)
	\end{dmath},
	\begin{dmath}
		\partial^n_\zeta \, \Phi(\zeta; \alpha, \beta, s) = {
		\uo^n \op{e}^{\zeta\uo} \,[\nu_{\alpha, \,\beta, s}] = 
		(\beta)_n \,\Phi(\zeta; \alpha + n, \beta + n, s)
	}
	\end{dmath}.
\end{dgroup}

We have already mentioned the importance of the Lerch transcendent function, descending from the fact that many special functions can be defined through it or can be derived as a particular case. For example, the inverse tangent integral, denoted by the symbol $\op{Ti}(\zeta)$ \cite{WeissteinInverseTangentIntegral},  is defined in terms of the generalised  Lerch function as:
\begin{dmath}
	2^s \frac{\op{Ti}(\zeta)	}{\zeta} = \Phi(-\zeta^2; \mbox{\small $\frac{1}{2}$}, 1, s)
\end{dmath}.
The umbral image of the r.h.s. of the previous identity is a Gaussian, namely:
\begin{dmath}\label{eq:umbral_ti}
	2^s \frac{\op{Ti}(\zeta)	}{\zeta} = \op{e}^{-\zeta^2 \uo}[\nu_{\frac{1}{2}, 1, s}]
\end{dmath}.
The successive derivatives with respect to $\zeta$ of both sides of \cref{eq:umbral_ti} can be easily obtained exploiting the following rule, valid for an ordinary Gaussian \cite{Babusci2019}:
\begin{dmath}\label{eq:gaussian_n_derivative}
	\partial^n_\zeta \op{e}^{a \zeta^2} = H_n(2 a \zeta, a)\op{e}^{a \zeta^2}
\end{dmath},
where $H_n(x,y)$ are the two-variable Hermite-Kamp{\'e} de F{\'e}ri{\'e}t polynomials: 
\begin{dmath}
	H_n(x, y) \defeq n! \sum_{r=0}^{\floor{\frac{n}{2}}}\frac{x^{n-2r}y^r}{(n-2r)!r!}
\end{dmath}.
Using \cref{eq:gaussian_n_derivative}, we eventually obtain:
\begin{dmath}\label{eq:umbral_ti}
	2^s \partial^n_\zeta \left(\frac{\op{Ti}(\zeta)}{\zeta}\right) = 
		H_n(-2 \zeta \uo, -\uo)\op{e}^{-\zeta^2 \uo}[\nu_{\frac{1}{2}, 1, s}] = {
		(-1)^n n! \sum_{r=0}^{\floor{\frac{n}{2}}}\frac{(-1)^r (2\zeta)^{n-2r}}{(n-2r)!r!}\uo^{n-r}\op{e}^{-\zeta^2 \uo}[\nu_{\frac{1}{2}, 1, s}] 
	} =
		(-1)^n n! \sum_{r=0}^{\floor{\frac{n}{2}}}\frac{(-1)^r (2\zeta)^{n-2r}}{(n-2r)!r!} (1)_{n-r} \,\Phi(-\zeta^2; n-r+1/2, 1+n-r, s)
\end{dmath}.
Before concluding this section, we would like to underline the link between the Lerch function and the polylogarithm function. Their entanglement is well known \cite{WeissteinInverseTangentIntegral} and the generalisations we have discussed so far offer an IUT view of the polylogarithm function, which writes:
\begin{dmath}
	\op{Li}_s(\zeta) = {
		\sum_{r=0}^\infty \frac{\zeta^r}{r^s} = 
		\zeta\,\Phi(\zeta; 1,1,s)
	}
\end{dmath}.
It is evident that
\begin{dmath}
	\frac{\op{Li}_s(\zeta)}{\zeta} = {
		\op{e}^{\zeta\uo} \,[\nu_{1, 1, s}]
	}
\end{dmath},
hence the properties of polylogarithms can be studied using this exponential $\nu_{1, 1, s}$-umbral image.
It is, for example, easily checked that
\begin{dmath}
	\frac{\op{Li}_s(\zeta_1 + \zeta_2)}{\zeta_1 + \zeta_2} = {
		\op{e}^{(\zeta_1 + \zeta_2)\uo} \,[\nu_{1, 1, s}] =
		\sum_{r=0}^\infty \frac{\zeta_1^r}{r!} \,r!\Phi(\zeta_2, 1+r,1+r,s)
	}
\end{dmath}.
Although we have so far only considered integer powers of the umbral operator $\uo$,  the expression
\begin{dmath}
	\uo^\lambda	[\nu_{\alpha, s}] = \frac{\Gamma(1+\lambda)}{(\lambda + \alpha)^s}
\end{dmath}
is well defined for any complex $\lambda$, provided $t = \lambda$ does not correspond to a singularity of the function $\nu_{\alpha, s}(t)$.
This enables to define "polylogarithms of non-integer order" such as:
\begin{dmath}
	g_s(\zeta) = {
		\frac{\op{Li}^{(1/2)}_s(\zeta)}{\zeta} =
		\op{e}^{\zeta\uo^{1/2}} \,[\nu_{1, 1, s}] = 
		\sum_{r=0}^\infty \frac{\Gamma(1 + r/2)}{(1+ r/2)^s}\,\frac{\zeta^r}{r!}
	}
\end{dmath}.
It is interesting to note that
\begin{dmath}\label{eq:identity_g_s}
	\int_{-\infty}^\infty \mathrm{d}x \op{e}^{-x^2} g_s(2 x \zeta) = \sqrt{\pi}\,\frac{\op{Li}_s(\zeta^2)}{\zeta^2}
\end{dmath}.
The above identity, whose proof is sketched in \cref{sec:appendix_1}, states that the polylogarithm function is the Gauss transform of its counterpart of order $1/2$.

We have underscored that that the defining series of the Lerch transcendent has a limited interval of convergence. We can however extend the domain of the function by analytic continuation. In particular, the following integral representation (see \cref{sec:appendix_2}):
\begin{dmath}\label{eq:lerch_extended}
	\Phi(\zeta, \alpha, s) = \frac{1}{\Gamma(s)}\int_{0}^\infty \mathrm{d}t \,\frac{t^{s - 1}x\op{e}^{-\alpha t}}{1 - \zeta \op{e}^{-t}} 
\end{dmath},
reproduces the series \cref{eq:lerch} for $|\zeta| \;< 1$, but converges for values of the variable exceeding unity, namely $\zeta \in \bbb{C} \backslash[1,\infty), \op{Re} s > 0, \op{Re} \alpha > 0$ (see Fig. (1)). The integral formula also holds in $\zeta = 1$ if $\op{Re} s > 1$.

\section{Legendre $\chi$ function}\label{sec:chi_legendre}

The $\chi$ function, introduced by Legendre in his book \emph{Exercice de Calcul Int{\'e}gral sur divers ordres de transcendantes et sur les quadratures}, has been the subject of influential research. A partial list of more recent studies is reported in refs. \cite{Lewin1958}. In modern notation, the $\chi$ function is defined as:
\begin{dmath}\label{eq:chi}
	\chi_s(\zeta) \defeq {
		\sum_{r=0}^\infty \frac{\zeta^{2r+1}}{(2r+1)^{s}} =
		\frac{\zeta}{2^s}\,\Phi(\zeta^2; \mbox{\small $\frac{1}{2}$}, s)
		\condition*{|\zeta|\; < 1}
	}
\end{dmath},
where we have made explicit its relationship with the Lerch transcendent.

In virtue of the following decomposition:
\begin{dmath}
	\Phi(\zeta; 1, s) = {
		\op{e}^{\zeta\uo} \,[\nu_{1, s}] =
		\cosh(\zeta\uo) \,[\nu_{1, s}] + \sinh(\zeta\uo) \,[\nu_{1, s}] =
		\op{c}_s(\zeta) + \op{s}_s(\zeta)
	}
\end{dmath},
where
\begin{dgroup}
	\begin{dmath}
		\op{c}_s(\zeta) \defeq {
			\cosh(\zeta\uo) \,[\nu_{1, s}] =
			\sum_{r=0}^\infty \frac{\zeta^{2r}}{(2r + 1)^s}
		}
	\end{dmath},
	\begin{dmath}
		\op{s}_s(\zeta) \defeq {
			\sinh(\zeta\uo) \,[\nu_{1, s}] =
			\sum_{r=0}^\infty \frac{\zeta^{2r+1}}{(2r+2)^s}
		}
	\end{dmath},
\end{dgroup}
we obtain by direct inspection the following $\nu_{1, s}$-umbral image for $\chi_s$:
\begin{dmath}
	\chi_s(\zeta) = {
		\zeta\,\op{c}_s(\zeta) = 
		\zeta\, \cosh(\zeta\uo) \,[\nu_{1, s}]
	}
\end{dmath}.

It is worth noting that:
\begin{dgroup}
	\begin{dmath}\label{eq:c_derivative}
		\partial_\zeta \op{c}_s(\zeta) = {
			\uo\,\sinh(\zeta\uo)\,[\nu_{1, s}] = 
			\sum_{r=0}^\infty \frac{2r+2}{(2r+3)^{s}}\,\zeta^{2r+1}
		}
	\end{dmath},
	\begin{dmath}\label{eq:s_derivative}
		\partial_\zeta \op{s}_s(\zeta) = {
			\uo\,\cosh(\zeta\uo)\,[\nu_{1, s}] = 
			\sum_{r=0}^\infty \frac{2r+1}{(2r+2)^{s}}\,\zeta^{2r}
		}
	\end{dmath}.
\end{dgroup}
Even though the use of the umbral formalism is not crucial for the derivation of the previous identities, it may simplify the following generalisations:
\begin{dgroup}
	\begin{dmath}\label{eq:c_n_derivative}
		\partial_\zeta^{2n} \op{c}_s(\zeta) = {
			\uo^{2n}\,\sinh(\zeta\uo)\,[\nu_{1, s}] = 
			\sum_{r=0}^\infty \frac{(2n+1)_{2r}}{(2r+2n +1)^{s}}\,\zeta^{2r}
		}
	\end{dmath},
	\begin{dmath}\label{eq:s_n_derivative}
		\partial_\zeta^{2n} \op{s}_s(\zeta) = {
			\uo^{2n}\,\cosh(\zeta\uo)\,[\nu_{1, s}] = 
			\sum_{r=0}^\infty \frac{(2n+1)_{2r+1}}{(2r+2n+2)^{s}}\,\zeta^{2r+1}
		}
	\end{dmath}.
\end{dgroup}

In analogy to the case of Lerch transcendent, the domain of definition of the Legendre chi function can be extended to a larger region using the following integral representation (see \cref{sec:appendix_2}):
\begin{dmath}\label{eq:chi_extended}
	\chi_s(\zeta) = \frac{\zeta}{\Gamma(s)}\int_{0}^\infty \mathrm{d}t \,\frac{t^{s - 1}\op{e}^{-t}}{1 - \zeta^2 \op{e}^{-2t}}
\end{dmath}.

\section{Polygamma function}\label{sec:polygamma}

The polygamma function of order $m$ ($m$-polygamma for short) is a meromorphic function on the complex plain,
defined as the $(m + 1)$-th derivative of the logarithm of the Gamma function:
\begin{dmath}\label{eq:polygamma}
	\psi^{(m)}(\alpha) \defeq {
		\partial_\alpha^m \,\psi(\alpha) = 
		\partial_\alpha^{m + 1} \,\ln \Gamma(\alpha) \condition*{m \in \bbb{N}^0}
	}
\end{dmath},
where $\psi^{(0)}(\alpha) = \psi(\alpha) = \Gamma'(\alpha)/\Gamma(\alpha)$ is the digamma function.
The $m$-polygamma is holomorphic on $\mathbb{C} \backslash \mathbb{Z} _{\leq 0}$, with poles of order $m+1$ at all the nonpositive integers.
Its relation with the Lerch transcendental is provided by the formula:
\begin{dmath}\label{eq:polygamma_lerch}
	\psi^{(m)}(\alpha) = (-1)^{m+1} m! \Phi(1; \alpha, m+1)
\end{dmath}.
It immediately follows from \cref{eq:umbral_lerch} that:
\begin{dmath}\label{eq:polygamma_lerch2}
	\psi^{(m)}(\alpha) = {
		(-1)^{m+1} m!\op{e}^\uo \, [\nu_{\alpha, m+1}] =
		(-1)^{m+1} m!\sum_{r=0}^\infty \frac{1}{(r+\alpha)^{m+1}} 
		} =
		(-1)^{m+1} m! \zeta(m+1, \alpha)
\end{dmath},
where $\zeta(s, \alpha) = \Phi(1; \alpha, s) = \op{e}^\uo \, [\nu_{\alpha, m+1}] $ is the Hurwitz zeta function.
By using in \cref{eq:polygamma_lerch} the identity
\begin{dmath}
	m! = {\Gamma(m+1) = \int_0^\infty \mathrm{d}x\, x^m \op{e}^{-x}}
\end{dmath},
after the change of integration variable $x \mapsto t(r+\alpha)$ and the interchange of sum and integral we obtain the integral representation:
\begin{dmath}\label{eq:polygamma_integral}
	\psi^{(m)}(\alpha) = {
		(-1)^{m+1}\int_0^\infty \mathrm{d}t\, \frac{t^m \op{e}^{-\alpha t}}{1-\op{e}^{-t}} 
		}
\end{dmath},
converging for $m > 0$ and $\op{Re} \alpha > 0$.

Further comments on the relevance and importance of the polygamma function for IUT will be discussed in the forthcoming conclusive section.

\section{Conclusions}\label{sec:conclusions}

One of the leitmotivs of the present investigation has been the possibility of extending the convergence of the series representative of the functions we have introduced, using different forms of the integral representation. This is a noticeable element of discussion because of its potential application in the context of summability.
As is well known the Euler series \cite{Cartier2000, Dattoli2010EulerLegacy}
\begin{dmath}\label{eq:euler_series}
	d_1(x) \defeq \sum_{r=0}^\infty (-1)^r r! x^r	
\end{dmath},
emerging from the perturbative solution of the equation
\begin{dmath}\label{eq:euler_eq}
	x^2 y' + y = x \condition*{y(0) = 0}	
\end{dmath},
has zero convergence radius.  However, using the integral representation of the factorial, the series in \cref{eq:euler_series} can be cast in the form of the integral representation, reported below:
\begin{dmath}\label{eq:euler_integral}
	d_1(x) = {\sum_{r=0}^\infty (-1x)^r \int_0^\infty \mathrm{d}t\,\op{e}^{-t} t^r	 \approx \int_0^\infty \mathrm{d}t\,\frac{\op{e}^{-t}}{1+xt} = \bar{d}_1(x)}
\end{dmath},
obtained after interchanging the symbols of summation and integral. This is an abuse (hence the symbol $\approx$, instead of the ordinary symbol of equality), since it holds for values of the variable $x$ allowing the convergence of the series in \cref{eq:euler_integral}. The consequence of this illegitimate procedure is twofold:
\begin{enumerate}
	\item $d_1 (x)$ has been associated with an integral transform converging for all positive $x$ values
	\item The solution of the differential equation in \cref{eq:euler_eq} can be written as
	\begin{dmath}
		y(x) = x \bar{d}_1(x)
	\end{dmath}.
\end{enumerate}
A further example of a similar manipulation is offered by the identities:
\begin{dmath}\label{eq:euler_series_2}
	d_2(x) \defeq {\sum_{r=0}^\infty (-1)^r (r!)^2 x^r \approx
		\int_0^\infty \mathrm{d}u\int_0^\infty \mathrm{d}v\,\frac{\op{e}^{-u+v}}{1+uvx} =
		\bar{d}_2(x)
	}
\end{dmath},
which state that an even more diverging series can be associated with a well-behaved function in the positive $x$ region. The price to be paid is the introduction of a double integral accounting for the squared factorial. The possible extension to any integer power in terms of multiple integral transforms is easily guessed.

A small step further is accomplished by considering the identities:
\begin{dmath}\label{eq:euler_series_2}
	d_2(x; \alpha, \beta) \defeq 
		{\sum_{r=0}^\infty (-1)^r (\Gamma(1+\beta +\alpha r)^2 x^r} 
		\approx {
		\int_0^\infty \mathrm{d}u\int_0^\infty \mathrm{d}v\,\frac{\op{e}^{-u+v}(uv)^{\,\beta}}{1+(uv)^\alpha x} = 
		\bar{d}_2(x; \alpha, \beta)
	}
\end{dmath},
suggesting the possibility of an extension of the IUT analysis for integer order Le Roy function, namely $L(x,\mu), \;\mu \in \bbb{N}^+$. 

The interest for this aspect of the problem is enhanced by the fact that $L(x,n)$ are Humbert type Bessel \cite{Humbert1930, DattoliEtAl} and by the possibility of introducing a multidimensional Borel-Le Roy transform as e. g.:
\begin{dgroup}
	\begin{dmath}
		\op{e}^{-x} = \int_0^\infty \mathrm{d}u\,\int_0^\infty \mathrm{d}v\,u^{\,\beta} v^{\delta}L(xu \op{e}^{\gamma v}; \alpha, \beta, \gamma, \delta)
	\end{dmath},
	\begin{dmath}
		L(x; \alpha, \beta, \gamma, \delta) = \sum_{r=0}^\infty \,\frac{x^r}{\Gamma(1+\beta + \alpha r)\Gamma(1+\delta + \gamma r)}
	\end{dmath}.
\end{dgroup}

Regarding the discussion associated to the Polygamma function, we should emphasize that its introduction seems to be extraneous to the umbral formalism outlined in this article. We did not use any umbra vacuum and image function to define their framing within the IUT framework. 
This can however be easily fixed, if we note that we can reinterpret the integral representation \cref{eq:polygamma_integral} as a special case of
\begin{dmath}\label{eq:2_var_polygamma_integral}
	\psi^{(m)}(\alpha, \zeta) = {
		-\partial_\alpha^m\int_0^\infty \mathrm{d}t\, \frac{\op{e}^{-\alpha t}}{1-\zeta\op{e}^{-t}} 
		}
\end{dmath},
which realizes a two variable Polygamma and reduces to the Lerch function. Indeed we find the obvious generalization of \cref{eq:polygamma_lerch}
\begin{dmath}\label{eq:2_var_polygamma_lerch}
	\psi^{(m)}(\alpha, \zeta) = \Gamma(m+1)\,\Phi(\zeta; \alpha, m+1)
\end{dmath}.
The translation of \cref{eq:2_var_polygamma_integral} in umbral terms is straightforward and will be discussed elsewhere, within a more general context.

These final comments yield an idea of the directions along which future researches can be developed and will be discussed by the present authors.

\appendix
\section{Appendix}\label{sec:appendix}
\renewcommand{\theequation}{\thesection.\arabic{equation}}
\setcounter{equation}{0}

\subsection{Derivation of \cref{eq:identity_g_s}}\label{sec:appendix_1}

In order to derive \cref{eq:identity_g_s}, we exploit the umbral identity
\begin{dmath}
	g_s(2x\zeta) = {
		\op{e}^{2x\zeta\uo^{1/2}} \,[\nu_{1, 1, s}]
	}
\end{dmath}
to write:
\begin{dmath}
	I(\zeta) = {
		\int_{-\infty}^\infty \mathrm{d}x \op{e}^{-x^2} g_s(2 x \zeta) = 
		\int_{-\infty}^\infty \mathrm{d}x \op{e}^{-x^2+2x\zeta\uo^{1/2}} \,[\nu_{1, 1, s}]
	}
\end{dmath}.
By "completing the square" we easily obtain
\begin{dmath}
	I(\zeta) = {
		\int_{-\infty}^\infty \mathrm{d}y \op{e}^{-y^2} \op{e}^{\zeta^2\uo}\,[\nu_{1, 1, s}] =
		\sqrt{\pi}\,\frac{\op{Li}_s(\zeta^2)}{\zeta^2}
	}
\end{dmath},
where $y \defeq x-\zeta\uo^{1/2}$ and use has been made of the umbral identity:
\begin{dmath}
	\op{e}^{\zeta^2\uo} \,[\nu_{1, 1, s}] = \frac{\op{Li}_s(\zeta^2)}{\zeta^2}
\end{dmath}.

\subsection{Derivation of \cref{eq:lerch_extended} and \cref{eq:chi_extended}}\label{sec:appendix_2}

In order to prove \cref{eq:lerch_extended}, we use the integral definition of the gamma function and write:
\begin{dmath}
	\Phi(\zeta, \alpha, s) \,\Gamma(s) = {
		\sum_{r=0}^\infty \frac{\zeta^r}{(r + \alpha)^s}\int_0^\infty \frac{\mathrm{d}x}{x}\,x^s \op{e}^{-x} 
	}	
\end{dmath}.
By changing the integration variable, $x \to t = x/(n+\alpha)$, and after interchanging the sum and integral, we obtain:
\begin{dmath}
	\Phi(\zeta, \alpha, s) \,\Gamma(s) = {
		\sum_{r=0}^\infty \int_0^\infty \frac{\mathrm{d}t}{t}\,t^s \zeta^r\op{e}^{-(r+\alpha)t} =
		\int_0^\infty \frac{\mathrm{d}t}{t}\,t^s \op{e}^{-\alpha t}\sum_{r=0}^\infty (\zeta \op{e}^{-t})^r 
	} =
		\int_{0}^\infty \mathrm{d}t \,\frac{t^{s - 1}\op{e}^{-\alpha t}}{1 - \zeta \op{e}^{-t}} 
\end{dmath}.

The derivation of \cref{eq:chi_extended} is completely analogous. 

\bibliographystyle{unsrt}
\bibliography{bibliography}


\end{document}